\documentclass[12pt,leqno]{article}
\usepackage[pdftex]{graphicx}
\usepackage{amsfonts, amsmath}

\newcounter{conjecture}
\newcounter{remark}
\newcounter{corollary}

\newenvironment{corollary}{\medskip{\bf Corollary \thecorollary.}
\addtocounter{corollary}{1}\em}{\rm}
\newtheorem{theorem}{Theorem}
\newtheorem{lemma}{Lemma}
\newtheorem{proposition}{Proposition}

\newcommand{\lll}{\label}
\newcommand {\rrr}[1]{(\ref{#1})}

\def \be{\begin{equation}}
\def \ee{\end{equation}}
\def \bt{\begin{theorem}}
\def \et{\end{theorem}}
\def \bc{\begin{corollary}}
\def \ec{\end{corollary}}
\def \bea{\begin{eqnarray}}
\def \eea{\end{eqnarray}}
\def \bas{\begin{eqnarray*}}
\def \eas{\end{eqnarray*}}

\def \noi{\noindent}
\def \aa{\alpha}

\def \pp{\phi}

\def \vski{\vspace{12pt}}
\def \hvski{\vspace{6pt}}
\def \ff{\infty}

\def \({\left(}
\def \){\right)}

\def \nn{\nonumber}

\def \bc{\begin{center} }
\def \ec{\end{center} }
\def \bs{\begin{slide} }
\def \es{\end{slide} }

\def\square{{\vcenter{\vbox{\hrule height.3pt
         \hbox{\vrule width.3pt height5pt \kern5pt
            \vrule width.3pt}
         \hrule height.3pt}}}}
\def\qed{{\hfill $\square$ \bigskip}}

\newcounter{cccases}
\newcommand{\ccases}[1]{\begingroup \refstepcounter{cccases} {\bf \fontsize{14}{16}\selectfont Case \thecccases }  \label{#1}\endgroup}

\begin{document}

\title{Simplified results on electric resistance on a distance-regular graph.}

\author{
\begin{tabular}{cc}
\textit{Jacobus Koolen} & \textit{Greg Markowsky} \\
jacobus\_koolen@yahoo.com & gmarkowsky@gmail.com \\
Department of Mathematics & Department of Mathematical Sciences \\
POSTECH & Monash University\\
Pohang, 790-784 & Victoria, 3800 \\
Republic of Korea & Australia
\end{tabular}}

\maketitle

\begin{abstract}
Simplifications of a result from \cite{markool} concerning the electric resistance between points in a distance-regular graph are given. In particular, we prove that the maximal resistance between points is bounded by twice the resistance between neighbors. We also show that if the constant is weakened, then a very simple proof can be given.
\end{abstract}

\section{Introduction}

In this note, we will consider the electric resistance between points in a graph; that is, we imagine that a graph $G$ is a circuit with each edge representing a wire with unit resistance, and the effective resistance measures the ease with which current moves between points (details will be presented in Section \ref{prelim}). In particular, we are interested in the following result, which was originally conjectured by Biggs in \cite{biggs2}:

\bt \lll{bigguy}
There is a universal constant $K$ such that if G is a distance-regular graph with degree at least 3 and diameter D then

\be
\max_j r_j = r_D \leq K r_1,
\ee
where $r_j$ is the electric resistance between any two vertices of distance $j$.
\et

This theorem shows that the class of distance-regular graphs have strong regularity properties with respect to the electric resistance metric. Biggs conjectured further as to the optimal value for $K$.

\begin{proposition} \label{1}
We may take $K=1+\frac{94}{101} \approx 1.931$ in Theorem \ref{bigguy}, and equality holds only for the Biggs-Smith graph.
\end{proposition}

The result was proved earlier in \cite{markool}, but the proof is rather long and technical, and relies heavily on a library of classification theorems for small distance-regular graphs. One of the purposes of this note is to give a much shorter and simpler proof of this result, using new techniques which were developed in \cite{marparkool} in order to prove the more difficult assertion that $K \searrow 1$ as the degree of the graph goes to $\ff$. The new proof is an improvement over the old, but still requires several classification results and a detailed analysis of several different cases. However, if we allow ourselves to accept the worse constant $K=3$, we will see that there is a very short and simple proof requiring no classification results whatsoever. As this proof is likely to be of more interest to non-specialists than the more difficult ones, and should even be accessible to those with no prior knowledge of distance-regular graphs, we will take the time to present it as well.

\begin{proposition} \label{2}
We may take $K=3$ in Theorem \ref{bigguy}.
\end{proposition}

In the next section we introduce the framework which was developed by Biggs in \cite{biggs2} which will allow us to prove these results. The ensuing section gives the proofs of Propositions \ref{1} and \ref{2}, and we close with a few remarks in the final section.

\section{Preliminaries} \label{prelim}

Let $G$ be a connected graph. To define the effective resistance $r_{u,v}$ between points $u, v$, we attach a battery of unit voltage between $u$ and $v$ and take the reciprocal of the current which flows through the graph when each edge is taken to have unit resistance. $r_{u,v}$ is then a metric on the graph, and this metric has a large number of important connections to random walks; \cite{doysne} is a highly elegant introduction to this concept. The distance $d(x,y)$ between any two vertices $x,y$ of $G$
is the length of a shortest path between $x$ and $y$ in $G$. The diameter of $G$ is the maximal distance
occurring in $G$ and we will denote this by $D = D(G)$.
For a vertex $x \in G$, define $K_i(x)$ to be the set of
vertices which are at distance $i$ from $x~(0\le i\le
D)$ where $D:=\max\{d(x,y)\mid x,y\in V(G)\}$ is the diameter
of $G$. In addition, define $K_{-1}(x):=\emptyset$ and $K_{D+1}(x)
:= \emptyset$. We write $x\sim_{G} y$ or simply $x\sim y$ if two vertices $x$ and $y$ are adjacent in $G$. A connected graph $G$ with diameter $D$ is called
{\em distance-regular} if there are integers $b_i,c_i$ $(0 \le i
\leq D)$ such that for any two vertices $x,y \in V(G)$ with $d(x,y)=i$, there are precisely $c_i$
neighbors of $y$ in
$K_{i-1}(x)$ and $b_i$ neighbors of $y$ in $K_{i+1}(x)$
(cf. \cite[p.126]{drgraphs}). In particular, any distance-regular graph $G$ is regular with valency
$k := b_0$ and we define $a_i:=k-b_i-c_i$ for notational convenience. Note that the definition implies that $c_{i+1} |K_{i+1}(x)| =
b_i |K_i(x)|$, so that in fact $|K_i(x)| = \frac{b_0 \ldots b_{i-1}}{c_1 \ldots c_i}$. A straightforward consequence of the definition is that (cf. \cite[Proposition 4.1.6]{drgraphs})\\

(i) $k=b_0> b_1\geq \cdots \geq b_{D-1}$;\\
(ii) $1=c_1\leq c_2\leq \cdots \leq c_{D}$;\\
(iii) $b_i\ge c_j$ \mbox{ if }$i+j\le D$.\\

Henceforth we work entirely with a distance-regular graph $G$ on $n$ vertices with associated intersection array $(b_0, b_1, \ldots, b_{D-1};c_1, \ldots, c_D)$, and we assume further that $k = b_0 \geq 3$. The {\it Biggs potentials} are defined recursively for $0 \leq i \leq D-1$ by

\bea \label{smile}
&& \pp_0=n-1
\\ \nn && \pp_i = \frac{c_i\pp_{i-1}-k}{b_i}
\eea

This recursive definition leads to the explicit value:

\be \label{nut}
\phi_i = k\Big(\frac{1}{c_{i+1}} + \frac{b_{i+1}}{c_{i+1}c_{i+2}} + \ldots + \frac{b_{i+1} \ldots b_{D-1}}{c_{i+1} \ldots c_{D}} \Big).
\ee

Note that

\begin{equation} \label{}
\begin{split}
\phi_{i-1} - \phi_i = k\Big(\Big(&\frac{1}{c_i}-\frac{1}{c_{i+1}}\Big) + \Big(\frac{b_{i}}{c_{i}c_{i+1}}-\frac{b_{i+1}}{c_{i+1}c_{i+2}}\Big) + \ldots \\ & + \Big(\frac{b_{i} \ldots b_{D-2}}{c_{i} \ldots c_{D-1}} - \frac{b_{i+1} \ldots b_{D-1}}{c_{i+1} \ldots c_{D}}\Big) + \frac{b_{i} \ldots b_{D-1}}{c_{i} \ldots c_{D}} \Big).
\end{split}
\end{equation}

Conditions (i) and (ii) above show that this quantity is positive, so $\phi_i$ is a strictly decreasing sequence. In \cite{biggs2}(or see \cite{markool}), the following was shown.

\begin{proposition} \label{3}
The resistance between two vertices of distance $j$ in $G$ is given by
\be \label{}
\frac{2\sum_{0\leq i<j}\pp_i}{nk}
\ee
\end{proposition}

Note that this proposition, together with Proposition \ref{1} and the fact that $\phi_0 = n-1$, show that

\begin{proposition} \label{}
The resistance $r_{u,v}$ between any vertices $u,v \in G$ satisfies

\begin{equation} \label{}
r_{u,v} < \frac{4}{k}.
\end{equation}
\end{proposition}

In light of Proposition \ref{3}, it is clear that Propositions \ref{1} and \ref{2} can be verified by proving that

\be \lll{prince}
\phi_1+ \ldots + \phi_{D-1} < K \phi_0.
\ee

This is what we will show.

\section{Proof of Propositions}

We begin by describing the general technique which will be used in both proofs. It is well-known that the case $b_1 = 1$ occurs only for cocktail party graphs, and the results are simple to verify in that case, so we will assume always that $b_1 \geq 2$. It is clear that \rrr{smile} implies

\begin{equation} \label{mover}
\pp_i < \frac{c_i}{b_i}\pp_{i-1}.
\end{equation}

This will be very useful to us so long as $b_i > c_i$. At such point as $c_i \geq b_i$, however, it will be more profitable to bound the expression in \rrr{nut}, as (iii) above implies that this occurs when $i$ is relatively close to $D$, and $\phi_i$ will therefore be the sum of a small number of small terms. In light of this, we set $j=\inf\{i:c_i \geq b_i\}$. We consider $\phi_0, \phi_1, \ldots, \phi_{j-1}$ to be the $head$ of the sequence, and $\phi_j, \ldots, \phi_{D}$ to be the $tail$. There is another interesting consequences of these definitions. Recall that $|K_i(x)| = \frac{b_0 \ldots b_{i-1}}{c_1 \ldots c_i}$, where $K_i(x)$The following lemma will be key for bounding the tail.

\begin{lemma} \label{run}
$\phi_j + \ldots + \phi_{D-1} \leq (j-1/2) \phi_{j-1}$.
\end{lemma}

This lemma is not particularly difficult, and a proof can be found in \cite{marparkool}. For bounding the head, we will simply observe that when $b_i > c_i$ then since $b_i \leq b_1$ and $b_i, c_i$ are integers we have $\frac{c_i}{b_i} < \frac{b_1-1}{b_1}$. In conjunction with \rrr{mover}, we then have

\begin{equation} \label{mover2}
\pp_i < \Big(\frac{b_1-1}{b_1}\Big)\pp_{i-1}.
\end{equation}

Furthermore, since $c_1 = 1$, \rrr{mover} implies

\begin{equation} \label{mover3}
\pp_1 < \frac{1}{b_1}\pp_{0}.
\end{equation}

We begin with the easier Proposition \ref{2}.

\vski

\noi {\bf Proof of Proposition \ref{2}:} To simplify the notation, we set $\aa = \frac{b_1-1}{b_1}$. We bound the head and tail in \rrr{prince} as described above, using \rrr{mover2}, \rrr{mover3}, and Lemma \ref{run}, to obtain

\begin{equation} \label{r}
\frac{\phi_1+ \ldots + \phi_{D-1}}{\phi_0} \leq \frac{1}{b_1} + \frac{\aa}{b_1} + \ldots + \frac{\aa^{j-2}}{b_1} + (j-1/2) \frac{\aa^{j-2}}{b_1}.
\end{equation}

Note that the term $\frac{\aa^{j-2}}{b_1}$ corresponds to $\frac{\phi_{j-1}}{\phi_0}$, so that the final term is the bound on $\frac{\phi_j+ \ldots + \phi_{D-1}}{\phi_0}$ given by Lemma \ref{run}. We now replace the head bound by the geometric series

\begin{equation} \label{h}
\frac{1}{b_1} + \frac{\aa}{b_1} + \frac{\aa^{2}}{b_1} + \ldots  = \frac{1}{b_1(1-\aa)} = \frac{1}{b_1(1-\frac{b_1-1}{b_1})} = 1.
\end{equation}
In order to control the tail term, we set $f(i) = \frac{(i-1/2)\aa^{i}}{b_1}$. Note that if $i \geq b_1$ then

\be \lll{}
\frac{f(i+1)}{f(i)} = \frac{b_1-1}{b_1} \times \frac{i+1/2}{i-1/2} \leq \frac{i-1}{i} \times \frac{i+1/2}{i-1/2} < 1,
\ee

whereas if $i \leq b_1-1$ we have

\be \lll{}
\frac{f(i+1)}{f(i)} = \frac{b_1-1}{b_1} \times \frac{i+1/2}{i-1/2} \geq \frac{i}{i+1} \times \frac{i+1/2}{i-1/2} > 1.
\ee

Thus, $f(i)$ attains its maximum at $i=b_1$. We therefore have

\begin{equation} \label{cp}
(j-1/2) \frac{\aa^{j-2}}{b_1} \leq (b_1-1/2) \frac{\aa^{b_1-2}}{b_1} \leq \frac{b_1-1/2}{b_1} < 1.
\end{equation}

Combining the estimates \rrr{r}, \rrr{h}, and \rrr{cp} gives the proposition. \qed

\vski

\noi {\bf Proof of Proposition \ref{1}:} The proof proceeds by considering a number of separate cases. We will show

\be
\frac{\phi_1+ \ldots + \phi_{D-1}}{\phi_0} < .93
\ee

\noi for all graphs other than the Biggs-Smith graph. If $b_1 \leq 2$, then it is known that either $D \leq 2$ or $k \leq 4$. We may therefore reduce our problem to the case $b_1 \geq 3$ by disposing of the following two cases.

\vski

\ccases{tit}: $D \leq 2$.

\hvski

There is nothing to show for $D=1$, and for $D=2$ we need only show $\phi_1<.93\phi_0$. This is clear from \rrr{mover3} and the assumption $b_1 \geq 2$. \qed

\ccases{deg3and4}: $k=3$ or $4$.

\hvski

The distance-regular graphs of valency 3 and 4 have been classified(see \cite[Thm 7.5.1]{drgraphs} and \cite{broukool}). The corresponding values of $\frac{\phi_1+ \ldots \phi_{D-1}}{\phi_0}$ are given in the following table, which contains only the graphs with $D \geq 3$. Note that all are less than $.93$ except for that of the Biggs-Smith graph.

\vski

\begin{tabular}{ l c c r }
Name & Vertices & Intersection array & $\frac{\phi_1+ \ldots \phi_{D-1}}{\phi_0}$ \\
\hline
Cube & 8 & (3,2,1;1,2,3)  &   0.428571 \\
Heawood graph & 14 & (3,2,2;1,1,3)  &  0.461538 \\
Pappus graph & 18 & (3,2,2,1;1,1,2,3)  &      0.588235 \\
Coxeter graph & 28 & (3,2,2,1;1,1,1,2)   &     0.666667 \\
Tutte's 8-cage & 30 & (3,2,2,2;1,1,1,3)    &    0.655172 \\
Dodecahedron & 20 & (3,2,1,1,1;1,1,1,2,3) &   0.842105 \\
Desargues graph & 20 & (3,2,2,1,1;1,1,2,2,3)  &  0.710526 \\
Tutte's 12-cage & 126 & (3,2,2,2,2,2;1,1,1,1,1,3)  &     0.872 \\
Biggs-Smith graph & 102 & (3,2,2,2,1,1,1;1,1,1,1,1,1,3) &  0.930693 \\
Foster graph & 90 & (3,2,2,2,2,1,1,1;1,1,1,1,2,2,2,3)  &      0.896067 \\
$K_{5,5}$ minus a matching& 10 & (4,3,1;1,3,4) &   0.296296\\
Nonincidence graph of $PG(2,2)$ &  14 & (4,3,2;1,2,4)  &  0.307692\\
Line graph of Petersen graph &  15 & (4,2,1;1,1,4)  &  0.428571\\
4-cube &  16 & (4,3,2,1;1,2,3,4) &       0.422222\\
Flag graph of $PG(2,2)$ & 21 & (4,2,2;1,1,2) &   0.5\\
Incidence graph of $PG(2,3)$ & 26 & (4,3,3;1,1,4)  &  0.32\\
Incidence graph of $AG(2,4)$-p.c. & 32 & (4,3,3,1;1,1,3,4)  &      0.376344\\
Odd graph $O_4$ & 35 & (4,3,3;1,1,2)  &  0.352941\\
Flag graph of $GQ(2,2)$ & 45 & (4,2,2,2;1,1,1,2)  &      0.681818\\
Doubled odd graph & 70& (4,3,3,2,2,1,1;1,1,2,2,3,3,4)&    0.521739\\
Incidence graph of $GQ(3,3)$ & 80 & (4,3,3,3;1,1,1,4) &       0.417722\\
Flag graph of $GH(2,2)$ & 189 & (4,2,2,2,2,2;1,1,1,1,1,2) &      0.882979\\
Incidence graph of $GH(3,3)$ & 728 & (4,3,3,3,3,3;1,1,1,1,1,4) &      0.485557\\
\end{tabular} \qed

\qed

\vski

\ccases{cm}: $b_1 \geq 3, c_2=1$.

\hvski

If $j = 2$, then applying Lemma \ref{run} gives $\phi_1 + \ldots + \phi_{D-1} \leq \frac{5}{2} \phi_1$. As $\phi_1 \leq \frac{\phi_0}{b_1} \leq \frac{\phi_0}{3}$, we get

\begin{equation} \label{}
\frac{\phi_1+ \ldots + \phi_{D-1}}{\phi_0} \leq \frac{5}{6} <.93
\end{equation}

If $j=3$, then $b_2 \geq 2$, hence $\phi_2 \leq \frac{\phi_1}{2} \leq \frac{\phi_0}{6}$. By Lemma \ref{run} we have

\begin{equation} \label{}
\frac{\phi_1+ \ldots + \phi_{D-1}}{\phi_0} \leq \frac{\phi_0/3 + (7/2)\phi_0/6}{\phi_0} = \frac{11}{12} < .93
\end{equation}

If $j \geq 4$ and $b_2=2$ then we must have $b_3=2, c_3=1$, so that $\frac{b_2 b_3}{c_2 c_3} = 4$. On the other hand, if this does not occur than $\frac{b_2}{c_2} \geq 3$. We will consider these cases separately.

\vski

{\bf Subcase 1:} $\frac{b_2}{c_2} \geq 3$.

\hvski

For $i < j$ we have $b_1 \geq b_i > c_i$, and for any $i$ with $c_i>1$ we must have $b_i<b_1$, by Proposition 5.4.4 in \cite{drgraphs}. Thus, $\frac{c_i}{b_i} \leq \frac{b_1-2}{b_1-1}$. Define $\aa = \frac{b_1-2}{b_1-1}$. Applying Lemma \ref{run} we have

\be \label{}
\frac{\phi_1+ \ldots + \phi_{D-1}}{\phi_0} \leq \frac{1}{b_1}+ \frac{1}{3b_1}+ \frac{\aa}{3b_1}+ \ldots + \frac{\aa^{j-3}}{3b_1} + \frac{(j-1/2)\aa^{j-3}}{3b_1}
\ee

Replace the second through $(j-1)$th term by a geometric series to obtain

\bea \label{kyl}
\frac{\phi_1+ \ldots + \phi_{D-1}}{\phi_0} < \frac{1}{b_1}+ \frac{1}{3b_1}\Big(\frac{1}{1-\frac{{b_1-2}}{b_1-1}}\Big) + \frac{(j-1/2)\aa^{j-3}}{3b_1}
\\ \nn = \frac{1}{b_1} + \frac{b_1-1}{3b_1} + \frac{(j-1/2)\aa^{j-3}}{3b_1}.
\eea

Using the same technique as in the proof of Proposition \ref{2}, set $f(i) = \frac{(i-1/2)\aa^{i-3}}{3b_1}$. Note that if $i \geq b_1-1$ then

\be \lll{}
\frac{f(i+1)}{f(i)} = \frac{b_1-2}{b_1-1} \times \frac{i+1/2}{i-1/2} \leq \frac{i-1}{i} \times \frac{i+1/2}{i-1/2} < 1,
\ee

whereas if $i \leq b_1-2$ we have

\be \lll{}
\frac{f(i+1)}{f(i)} = \frac{b_1-2}{b_1-1} \times \frac{i+1/2}{i-1/2} \geq \frac{i}{i+1} \times \frac{i+1/2}{i-1/2} > 1
\ee

Thus, $f(i)$ attains its maximum at $i=b_1 - 1$. Using this to bound the final term in \rrr{kyl} for $b_1 \geq 4$ gives

\bea \label{}
\frac{\phi_1+ \ldots + \phi_{D-1}}{\phi_0} < \frac{1}{b_1} + \frac{b_1-1}{3b_1} + \frac{(b_1-3/2)\aa^{b_1-4}}{3b_1}
\\ \nn \leq  \frac{1}{b_1} + \frac{b_1-1}{3b_1} + \frac{1}{3} = \frac{2b_1+2}{3b_1} < .93
\eea

If $b_1=3$, then since $j \geq 4$ we can simply plug in $j=4$ to get

\bea \label{}
\frac{\phi_1+ \ldots + \phi_{D-1}}{\phi_0} < \frac{1}{b_1} + \frac{b_1-1}{3b_1} + \frac{(7/2)\aa}{3b_1}
\\ \nn = \frac{3}{4} < .93
\eea

{\bf Subcase 2:} $\frac{b_2 b_3}{c_2 c_3} \geq 4$.

\hvski

This follows much as in the previous case. Let $\aa = \frac{b_1-2}{b_1-1}$. Since $b_2 \geq b_3$ and $c_2 \leq c_3$ we must have $\frac{b_2}{c_2} \geq 2$. We then have

\be \label{jb}
\frac{\phi_1+ \ldots + \phi_{D-1}}{\phi_0} \leq \frac{1}{b_1}\!+ \! \frac{1}{2b_1}\!+ \! \frac{1}{4b_1}\!+ \! \frac{\aa}{4b_1}\!+ \! \ldots \!+ \! \frac{\aa^{j-4}}{4b_1}\! + \!\frac{(j-1/2)\aa^{j-4}}{4b_1}.
\ee

If $j=4$, then in fact the terms containing $\aa$'s are not present, and we get

\be \label{}
\frac{\phi_1+ \ldots + \phi_{D-1}}{\phi_0} \leq \frac{1}{b_1} + \frac{1}{2b_1} + \frac{1}{4b_1} + \!\frac{(4-1/2)}{4b_1} = \frac{10 \frac{1}{2}}{4b_1} \leq \frac{10 \frac{1}{2}}{12} < .93.
\ee

If $j>4$, then we again set $f(i) = \frac{(i-1/2)\aa^{i-4}}{4b_1}$, and use the argument from the previous subcase to conclude that $f(i)$ is decreasing for $i \geq b_1 - 1$ but increasing for $i< b_1-1$. If $b_1 < 6$, we may therefore replace $j$ by 5 in \rrr{jb} and sum the geometric series in $\aa$ to get a bound of

\be \label{}
\frac{\phi_1+ \ldots + \phi_{D-1}}{\phi_0} \leq \frac{3}{2b_1} + \frac{1}{4}\Big(\frac{b_1-1}{b_1}\Big) + \!\frac{(4 \frac{1}{2})\aa }{4b_1}.
\ee

If $b_1 = 4,5$, then this expression is seen to be less than $.93$ upon replacing $\aa$ by 1, while for $b_1 = 3$ a sufficient bound is obtained by using $\aa = \frac{1}{2}$. If $b_1 \geq 6$, we can replace $j$ by $b_1-1$ in \rrr{jb} and again sum the geometric series to obtain

\be \label{}
\frac{\phi_1+ \ldots + \phi_{D-1}}{\phi_0} \leq \frac{3}{2b_1} + \frac{1}{4}\Big(\frac{b_1-1}{b_1}\Big) + \!\frac{(b_1 - 3/2)\aa^{b_1 - 5} }{4b_1} < \frac{3}{2b_1} + \frac{1}{2} < .93.
\ee

\qed

\ccases{jr}: $b_1 \geq 3, j=3, c_2>1$.

\hvski

By Theorem 5.4.1 in \cite{drgraphs}, $c_2 \leq \frac{2}{3} c_3$. Suppose first that $c_3 > b_3$; then $D \leq 2j-1 = 5$ by Property (iii) in Section \ref{prelim}, and $\phi_2 < \frac{b_1-1}{b_1} \phi_1 < \frac{b_1-1}{b_1^2} \phi_0$. We therefore have

\begin{equation} \label{}
\frac{\phi_1+ \ldots + \phi_{D-1}}{\phi_0} < \frac{1}{b_1} + 3 \Big(\frac{b_1-1}{b_1^2}\Big) = \frac{4b_1-3}{b_1^2}.
\end{equation}

For $b_1 \geq 4$, this is less than $.93$. If $b_1=3$, then we must have $c_2=1$, since if $c_2>1$ then $b_2 < 3$ by Proposition 5.4.4 in \cite{drgraphs}, contradicting $j=3$. Thus, $\phi_2 < \frac{1}{2} \phi_1 < \frac{1}{6} \phi_0$, and we obtain

\begin{equation} \label{}
\frac{\phi_1+ \ldots + \phi_{D-1}}{\phi_0} < \frac{1}{3} + 3 \times \frac{1}{6} < .93.
\end{equation}

If $c_3=b_3 \leq b_2$, then if we assume $\frac{c_2}{b_2} \leq \frac{1}{2}$ by Lemma \ref{run} we have

\be \label{}
\frac{\phi_1+ \ldots + \phi_{D-1}}{\phi_0} \leq \frac{\phi_1+ (7/2)\phi_{2}}{\phi_0} \leq \frac{1}{b_1}+ \frac{7}{4b_1} = \frac{11}{4b_1} < .93
\ee

On the other hand, if it is not the case that $\frac{c_2}{b_2} \leq \frac{1}{2}$, then the proof of Theorem 5.4.1 of \cite{drgraphs} implies that $G$ contains a quadrangle. By Corollary 5.2.2 in \cite{drgraphs}, $D\leq\frac{2k}{k+1-b_1}$. It is straightforward to verify that the fact that $k \geq b_1+1$ implies that

\be \label{}
\frac{2k}{k+1-b_1} \leq b_1+1
\ee

We therefore see that the fact that $G$ contains a quadrangle implies $D \leq b_1+1$. Furthermore, we still have $\frac{c_2}{b_2} \leq \frac{2}{3}$ by Theorem
5.4.1 of \cite{drgraphs}. We therefore have

\be \label{}
\frac{\phi_1+ \ldots + \phi_{D-1}}{\phi_0} \leq \frac{\phi_1+ (b_1-1)\phi_{2}}{\phi_0} \leq \frac{1}{b_1}+ \frac{2(b_1-1)}{3b_1} = \frac{2b_1+1}{3b_1} < .93
\ee \qed

\ccases{the}: $b_1 \geq 3, j \geq 4, c_2 > 1, G$ contains a quadrangle.

\hvski

As in the argument given in Case \ref{jr}, we see that $G$ containing a quadrangle implies $D \leq b_1+1$. Furthermore,
Theorem 5.4.1 of \cite{drgraphs} implies that $c_3 \geq (3/2)c_2$. Since $j \geq 4$ and thus $b_2 \geq b_3 > c_3$
we must have $\frac{c_2}{b_2} \leq \frac{2}{3}$. This gives

\be \label{}
\frac{\phi_1+ \ldots + \phi_{D-1}}{\phi_0} \leq \frac{1}{b_1}+ (b_1-1)\frac{2}{3b_1} = \frac{2b_1+1}{3b_1} <.93
\ee

\qed

\ccases{guy}: $b_1 \geq 3, j \geq 4, c_2 > 1, G$ does not contain a quadrangle.

\hvski

In this case $G$ is a Terwilliger graph. By Corollary 1.16.6 of \cite{drgraphs}, if $k<50(c_2-1)$ then $D \leq 4$, which was covered in \cite{biggs2}. Thus, we can assume $k \geq 50(c_2-1) \geq 50$, which implies $b_1 \geq \frac{k}{2}>20$. By Theorem 3 of \cite{marparkool},

\begin{equation} \label{}
\phi_2 + \ldots + \phi_{D-1} < 9 \phi_1
\end{equation}

We then have

\begin{equation} \label{}
\frac{\phi_1+ \ldots + \phi_{D-1}}{\phi_0} < \frac{10 \phi_1}{\phi_0} < \frac{10 (\frac{\phi_0}{20})}{\phi_0} = \frac{1}{2}
\end{equation}
\qed

\bibliographystyle{alpha}
\bibliography{drgraphs}

\end{document}